\documentclass[a4paper,12pt,leqno]{amsart}
\usepackage{amsmath}
\usepackage{amsfonts}
\usepackage{amssymb}
\usepackage{mathrsfs}

\newcommand{\R}{\mathbb{R}}

\newcommand{\N}{\mathbb{N}}
\newcommand{\spt}{\operatorname{spt}}

\newcommand{\diam}{\operatorname{diam}}
\newcommand{\guio}[1]{\nobreakdash-\hspace{0pt}#1}

\swapnumbers \theoremstyle{plain}
\newtheorem{thm}[equation]{Theorem}

\theoremstyle{definition}

\theoremstyle{remark}
\newtheorem*{gracies}{Acknowledgements}

\numberwithin{equation}{section}

\title{Convergence of singular integrals with general measures}
\author{Pertti Mattila and Joan Verdera}
\date{}

\begin{document}

\begin{abstract}
We show that $L^2$-bounded singular integrals in metric spaces with
respect to general measures and kernels converge weakly. This
implies a kind of average convergence almost everywhere. For
measures with zero density we prove the almost everywhere existence
of principal values.
\end{abstract}

\maketitle

\section{Introduction}
Singular integrals with respect to general measures in $\R^n$, and
also in metric spaces, have been studied widely, see, e.g.,
\cite{C}, \cite{CW}, \cite{D1}, \cite{DS}, \cite{M}, \cite{P},
\cite{Ve} and \cite{V}. In this paper our setting is a separable
metric space $(X,d)$ with a finite Borel measure $\mu$ and a Borel
measurable antisymmetric kernel $K\colon X \times X \setminus
\{(x,y):x=y\}\rightarrow \R$. Antisymmetry means that
$$
K(x,y)=-K(y,x) \text{ for } x,y\in X,\quad x\not=y.
$$
Moreover, we shall assume that $K$ is bounded in $\{(x,y)\in X\times
X:d(x,y)>\delta\}$ for every $\delta>0$. We shall also always assume
that Vitali's covering theorem is valid for $\mu$ and the family of
closed balls. Although this is not automatically true even when $X$
is compact, it is true for example if $X=\R^n$ or $\mu$ is doubling,
see, e.g., \cite[Section 2.8]{F}.

The singular integral operator $T$ associated with $\mu$ and $K$ is
formally given by
$$
T(f)(x) = \int K(x,y)f(y)\,d\mu y.
$$
The problem which appears already in all classical cases such as the
Hilbert transform on $\R$, i.e., $K(x,y) = 1/(y-x)$, is that usually
this integral does not exist when $x \in \spt\mu$, the support of
$\mu$. When $\mu$ is the Lebesgue measure $\mathcal{L}^n$ on $\R^n$
and $K$ is a standard Calder\'on-Zygmund kernel, this can be
overcome by defining
\begin{equation}\label{eq:1.1}
T(f)(x) = \lim_{\epsilon \rightarrow 0} T_{\epsilon}f(x),
\end{equation}
where
$$
T_{\epsilon}(f)(x) = \int_{X \setminus
B(x,\epsilon)}K(x,y)f(y)\,d\mu y.
$$
Here $B(x,\epsilon)$ is the open ball with centre $x$ and radius
$\epsilon$. In such a case the limit exists trivially for smooth
functions due to cancellations, and by the denseness of smooth
functions in $L^1(\mathcal{L}^n)$ standard techniques can be used to
show that it exists almost everywhere for $L^1$-functions~$f$. For
general measures this approach fails. Unless $\mu$~has strong
symmetry properties around points in its support there are not
enough cancellations to guarantee the existence of the limit even
for constant functions. However, when $K$ is antisymmetric one often
defines $T(f)$ as a distribution by
\begin{equation}\label{eq:1.2}
(T(f),g) = (1/2)\iint K(x,y)(f(x)g(y) - f(y)g(x))\,d\mu x\,d\mu y
\end{equation}
when $f$ and $g$ are bounded Lipschitz functions, see \cite{C} or
\cite{D1}.

A central concept in the theory of singular integrals is the
boundedness in $L^2$. This can be formulated in several ways which
all agree in the classical case of Calder\'on-Zygmund kernels and
the Lebesgue measure. One way is to say that the distributionally
defined operator $T$, as in \eqref{eq:1.2}, is bounded in $L^2(\mu)$
if it has a bounded extension to~$L^2(\mu) \rightarrow L^2(\mu)$.
Another way is to require that the truncated operators
$T_{\epsilon}$, $\epsilon>0$, are uniformly bounded in $L^2(\mu)$.
This agrees very generally with the boundedness in $L^2(\mu)$ of the
sublinear maximal operator $T^{\ast}$:
\begin{equation}\label{eq:1.3}
T^{\ast}(f)(x) = \sup_{\epsilon>0}|T_{\epsilon}(f)(x)|,
\end{equation}
see \cite{NTV}.

A natural question is whether the $L^2$-boundedness forces the limit
$\lim_{\epsilon \rightarrow 0}T_{\epsilon}(f)(x)$ to exist for $\mu$
almost all $x\in X$. One would expect this to be true at least if
$\mu$ is an $m$-dimensional Ahlfors-David-regular measure in $\R^n$:
$$
r^m/C \leq \mu(B(x,r)) \leq Cr^m \text{ for } x \in \spt\mu,\quad
0<r<\diam(\spt\mu),
$$
and $K$ is the vector-valued Riesz kernel $|x-y|^{-m-1}(x-y)$. In
fact, by a result of Tolsa, see \cite{T1}, this is true when $m=1$
even for much more general measures, but the proof is based on very
special relations with the kernel $x/|x|^2$ (essentially the Cauchy
kernel $1/z$ for $z \in \mathbb{C} = \R^2$) and the so-called Menger
curvature. We shall discuss some relations of this problem to
rectifiability at the end of the paper. And we shall mention some
kernels for which $L^2$-boundedness does not give the almost
everywhere convergence of principal values.

In this paper we prove some substitutes for \eqref{eq:1.1} under the
$L^2$\guio{boundedness}:

\begin{thm}\label{thm:1.4}
Suppose that $T^{\ast}$ (defined by \eqref{eq:1.3}) is bounded in
$L^2(\mu)$, that is, there exists a constant $C_0$ such that
\begin{equation}\label{eq:1.5}
\int T^{\ast}(f)^2d\mu \leq C_0\int f^2\,d\mu
\end{equation}
for $f \in L^2(\mu)$. Then the truncated operators $T_{\epsilon}$
converge weakly in $L^2(\mu)$, that is, there exists a bounded
linear operator $T\colon L^2(\mu)\to L^2(\mu)$ such that
$$
\lim_{\epsilon\to 0}\int T_{\epsilon}(f)g\,d\mu=\int T(f)g\,d\mu
$$
for $f,g \in L^2(\mu)$. Moreover,
$$
T(f)(z) = \lim_{r \rightarrow 0} \frac{1}{\mu(B(z,r))} \int_{B(z,r)}
\left(\int_{X \setminus B(z,r)}K(x,y)f(y)\,d\mu y \right) \,d\mu x
$$
for $\mu$ almost all $z \in X$.
\end{thm}

So even if we don't know that $T(f)$ would exist as the limit of the
simpler integrals $T_{\epsilon}(f)$, we know that it is almost
everywhere the limit of the more complicated but still concrete
integrals of Theorem~\ref{thm:1.4}.

Observe that with some natural estimates the limit operator $T$
satisfies \eqref{eq:1.2}. This is so if, for example,
$$
\iint|K(x,y|\,d(x,y)\,d\mu y\,d\mu x < \infty,
$$
as one easily checks. In many cases also the converse in the first
part of Theorem~\ref{thm:1.4} is true. Namely, by the
Banach-Steinhaus theorem the weak convergence implies that the
truncated operators $T_{\epsilon}$ are uniformly bounded and, as
said before, often this is equivalent to the $L^2$-boundedness of
$T^{\ast}$.

We prove Theorem~\ref{thm:1.4} in Section~2. We first establish the
weak convergence. Then we deduce from it the average convergence
using Lebesgue differentiation theorem. We shall also indicate in
Section~3 another way of getting the average convergence via the
martingale convergence theorem.

In Section~4 we apply Theorem~\ref{thm:1.4} to prove the following
result on the existence of principal values for measures with zero
density:

\begin{thm}\label{thm:1.6}
Suppose $X=\R^n$ or $\mu$ is doubling. Let $h\colon
(0,\infty)\to(0,\infty)$ be an increasing function such that
$\lim_{r\to0}h(r)=0$, $h(2r)\leq Ch(r)$ for $r>0$ and that for
$x,y\in X$, $x\not=y$,
\begin{equation}\label{eq:1.7}
|K(x,y)|\leq\frac{1}{h(d(x,y))},
\end{equation}
and for $z\in X$, $z\not=x$ with $d(x,y)>2d(y,z)$,
\begin{equation}\label{eq:1.8}
|K(x,y)-K(x,z)|\leq\frac{d(y,z)}{d(x,y)h(d(x,y))}.
\end{equation}
Suppose also that for all $x\in X$ and $r>0$,
\begin{equation}\label{eq:1}
\mu(B(x,r))\leq h(r)
\end{equation}
and for $\mu$ almost all $x\in X$,
\begin{equation}\label{eq:1.9}
\lim_{r\to0}\frac{\mu(B(x,r))}{h(r)}=0.
\end{equation}
If $T^{\ast}\colon L^2(\mu)\to L^2(\mu)$ is bounded, then for $f\in
L^1(\mu)$ and for $\mu$ almost all  $x\in X$,
$$
\lim_{\epsilon\to0}T_{\epsilon}(f)(x)=T(f)(x)
$$
where $T$ is the weak limit operator of Theorem~~\ref{thm:1.4}.
\end{thm}

Note that originally $T(f)$ was only defined for $f\in L^2(\mu)$,
but under the assumptions of the theorem it has a unique extension
to $L^1(\mu)$ because we have the weak $L^1$-inequality: for $t>0$,
\begin{equation}\label{eq:1.10}
\mu(\{x\in X: |T^{\ast}(f)(x)|>t\})\leq C||f||_1/t.
\end{equation}
For the doubling measures in metric spaces this was proved in
\cite{CW} and for general measures in $\R^n$ in \cite{NTV}. The
assumptions on the kernels in \cite{NTV} are not quite same as above
but it is easy to check that the proofs can be modified.

Rather often the growth condition \eqref{eq:1} is a consequence of
the $L^2$ boundedness of $T^{\ast}$ (see  \cite [p.~56]{D1}).

For general kernels $K$ as above the assumption \eqref{eq:1.9} is
necessary as an example of David, which we discuss at the end of the
paper, shows.

A particular but interesting instance of the above result arises in
the following situation. We take $X= \R^n$ and an underlying measure
$\mu$ which satisfies the growth condition $ \mu(B(x,r))\leq
C\,r^m\,,$ for each $x$ and each $r > 0$\,. The kernel is a standard
smooth antisymmetric m-dimensional kernel satisfying the usual
conditions
\begin{equation*}
|K(x,y)|\leq \frac{1}{|x-y|^m},
\end{equation*}
and
\begin{equation}
|K(x,y)-K(x,z)|\leq \frac{|y-z|}{|x-y|^{m+1}}\,,\quad  |x-y|>2 |y-z|
\,.
\end{equation}
Then \eqref{eq:1.9} says that m-dimensional density vanishes for
$\mu$ almost all $x$, namely,
$$\lim_{r\rightarrow 0} \frac{\mu(B(x,r))}{r^m} = 0\,.$$
This, of course, excludes m-dimensional Ahlfors-David regular sets.
See the remarks in section 5\,.

\section{Proof of Theorem 1.4}
Let $B$ be a closed ball in $X$. We denote by $\chi_A$ the
characteristic function of a set $A$ and by $A^c$ its complement in
$X$. We have for all~$\epsilon>0$ ($1$ denotes the constant function
identically $1$),
$$
\int T_{\epsilon}(1) \chi_B \,d\mu = - \int T_{\epsilon}(\chi_B)
\,d\mu = - \int _ {B^c} T_{\epsilon}(\chi_B) \,d\mu,
$$
because by antisymmetry
$$
\int _ {B} T_{\epsilon}(\chi_B) \,d\mu = 0.
$$

Clearly, for all $x\in B^c$ there is the limit (since $B$ is closed)
$$
T(\chi_B)(x):=\lim_{\epsilon\to 0}T_{\epsilon}(\chi_B)(x).
$$

As $|T_{\epsilon}(\chi_B)|\leq T^{\ast}(\chi_B)\in L^1(\mu)$, the
dominated convergence theorem yields that
\begin{equation}\label{eq:2.1}
\lim_{\epsilon\to 0}\int T_{\epsilon}(1)\chi_B\,d\mu=
-\lim_{\epsilon\to0}\int_{B^c}T_{\epsilon}(\chi_B)\,d\mu=-\int_{B^c}T(\chi_B)\,d\mu.
\end{equation}

Call $S$ the dense subspace of $L^2(\mu)$ consisting of finite
linear combinations of characteristic functions of closed balls. (It
is easy to verify that $S$ is dense since we assumed Vitali's
covering theorem for $\mu$.) Fix~$f$ in $L^2(\mu)$ and take $b$ in
$S$ extremely close to $f$ in $L^2(\mu)$. Then for
$0<\epsilon<\delta$,
$$
\int (T_{\delta}(1)-T_{\epsilon}(1)) f\, d\mu\! =\! \int
(T_{\delta}(1)-T_{\epsilon}(1)) (f-b) \,d\mu + \int
(T_{\delta}(1)-T_{\epsilon}(1)) b \,d\mu .
$$

By \eqref{eq:2.1}, the second term goes to $0$ as $\delta\to0$. For
the first term we have by the Schwartz inequality and the
$L^2$-boundedness \eqref{eq:1.5} of $T^{\ast}$,
\begin{multline*}
\left|\int (T_{\delta}(1)-T_{\epsilon}(1))(f-b)\,d\mu\right| \leq
||T_{\delta}(1)-T_{\epsilon}(1)||_2||f-b||_2\\
\leq2||T^{\ast}(1)||_2||f-b||_2 \leq
2(C_0\mu(X))^{\frac{1}{2}}||f-b||_2 ,
\end{multline*}
which we can make as small as we want. This gives that the finite
limit
$$
\lim_{\epsilon\to 0}\int T_{\epsilon}(1)f \,d\mu
$$
exists for all $f\in L^2(\mu)$.

Let again $B$ be a closed ball and $f\in L^2(\mu)$. Then for
$\epsilon>0$,
\begin{multline*}
\int T_{\epsilon}(\chi_B)f\,d\mu =
\int_B\int_{B\setminus B(x,\epsilon)}K(x,y)\,d\mu yf(x)\,d\mu x\\
+\int_{B^c}\int_{B\setminus B(x,\epsilon)}K(x,y)\,d\mu yf(x)\,d\mu
x.
\end{multline*}
Applying what we proved above to the measure $\chi_B\mu$ we conclude
that the first integral converges as $\epsilon\to0$. The second
integral converges again by the dominated convergence theorem, since
$$
\left|\int_{B\setminus B(x,\epsilon)}K(x,y)\,d\mu yf(x)\right|\leq
T^{\ast}(\chi_B)(x)|f(x)|
$$
and $T^{\ast}(\chi_B)f\in L^1(\mu)$. Then also
$$
\lim_{\epsilon\to0}\int T_{\epsilon}(b)f\,d\mu
$$
exists for all $f\in L^2(\mu), b\in S$. Arguing as above with the
$L^2$\guio{boundedness} we find that
$$
\lim_{\epsilon\to0}\int T_{\epsilon}(g)f\,d\mu
$$
exists for all $f,g\in L^2(\mu)$. This yields easily that there
exists a bounded linear operator $T\colon L^2(\mu)\to L^2(\mu)$ such
that
$$
\int T(g)f\,d\mu=\lim_{\epsilon\to0}\int T_{\epsilon}(g)f\,d\mu
$$
for all $f,g\in L^2(\mu)$, and we have established the required weak
convergence.

Let $B=B(z,r)$ be an open ball with $\mu(B)>0$. Using the
antisymmetry of $K$ we have for all $\epsilon>0$,
\begin {equation*}
\begin{split}
\int_B T_{\epsilon}(f\chi_{B^c})\,d\mu &= - \int f \chi_{B^c}
T_{\epsilon}(\chi_B) \,d\mu\\*[5pt] &=-\int
fT_{\epsilon}(\chi_B)\,d\mu + \int
f\chi_BT_{\epsilon}(\chi_B)\,d\mu\\*[5pt]
&=\int_BT_{\epsilon}(f)\,d\mu +
\int_B(f-f_B)T_{\epsilon}(\chi_B)\,d\mu,
\end{split}
\end{equation*}
where $f_B=\frac{1}{\mu(B)}\int_Bf\,d\mu$ and
$\int_BT_{\epsilon}(\chi_B)\,d\mu=0$. Letting $\epsilon\to0$,  we
obtain for the weak limit operator $T$,
$$
\int_BT(f\chi_{B^c})\,d\mu = \int_BT(f)\,d\mu +
\int_B(f-f_B)T(\chi_B)\,d\mu.
$$
Dividing with $\mu(B)=\mu(B(z,r))$ and letting $r\to0$, we have for
$\mu$ almost all $z$ for the first term of the right hand side by
the Lebesgue differentiation theorem,
$$
\lim_{r\to0}\frac{1}{\mu(B(z,r))}\int_{B(z,r)}T(f)\,d\mu = T(f)(z),
$$
and for the second term by the Schwartz inequality,
$L^2$-boundedness of~$T$ and the Lebesgue differentiation theorem,
$$
\lim_{r\to0}\frac{1}{\mu(B(z,r))}\int_{B(z,r)}(f-f_{B(z,r)})T(\chi_{B(z,r)})\,d\mu
= 0.
$$
On the other hand,
$$
T_{\epsilon}(f\chi_{B^c})(x) = \int_{B^c\setminus
B(x,\epsilon)}K(x,y)f(y)\,d\mu y \to \int_{B^c}K(x,y)f(y)\,d\mu y
$$
as $\epsilon\to0$ for $x\in B$ with $|T_{\epsilon}(f\chi_{B^c})(x)|
\leq |T^{\ast}(f\chi_{B^c})(x)|$, and so by the dominated
convergence theorem,
$$
\int_BT(f\chi_{B^c})\,d\mu =
\lim_{\epsilon\to0}\int_BT_{\epsilon}(f\chi_{B^c})\,d\mu =
\int_B\int_{B^c}K(x,y)f(y)\,d\mu y\,d\mu x.
$$
Combining the above equations, we obtain
$$
\lim_{r\to0}\frac{1}{\mu(B(z,r))}\int_{B(z,r)}\int_{B(z,r)^c}K(x,y)f(y)\,d\mu
y\,d\mu x = Tf(z)
$$
for $\mu$ almost all $z\in X$. This proves the theorem.

For further reference we record for every ball $B$,
\begin{equation}\label{eq:2.2}
\int_BT(1)\,d\mu =\int_BT(\chi_{B^c})\,d\mu
=-\int_{B^c}T(\chi_B)\,d\mu
\end{equation}
which follows by antisymmetry.

\section{Martingales}

We introduce a general nested system of sets. Standard examples are
dyadic lattices of cubes in $\R^n$. For each $k \in
\N=\{1,2,\dotsc\}$ let $\mathcal{D}_k$~be a countable disjoint
partition of $X$ into $\mu$ measurable sets $D$ such that
$\mu(\partial D) = 0$. Let $\mathcal{D} =
\cup_{k=1}^\infty\mathcal{D}_k$. We assume that the
system~$\{\mathcal{D}_k\}$ is nested in the sense that every $D \in
\mathcal{D}_{k+1}$ is contained in some $D' \in \mathcal{D}_k$. Then
every $D' \in \mathcal{D}_k$ is a disjoint union of sets in
$\mathcal{D}_{k+1}$.

Suppose that $T^{\ast}$ is bounded in $L^2(\mu)$. Let $f\in
L^2(\mu)$ and $D \in \mathcal{D}_k$. As $\mu(\partial D) = 0$ we
have for $\mu$ almost all $x\in D$,
$$
\int_{D^c}K(x,y)f(y)\,d\mu y = \lim_{\epsilon\to0}\int_{D^c\setminus
B(x,\epsilon)}K(x,y)f(y)\,d\mu y.
$$
Moreover,
$$
\left|\int_{D^c\setminus B(x,\epsilon)}K(x,y)f(y)\,d\mu y\right|\leq
T^{\ast}(f\chi_{D^c})(x)\leq T^{\ast}(f)(x)+T^{\ast}(f\chi_{D})(x).
$$
If also $g\in L^2(\mu)$ we get by the dominated convergence theorem
\begin{multline}\label{eq:3.1}
\int_D\left|\int_{D^c}K(x,y)f(y)\,d\mu yg(x)\right|\,d\mu x\\
\leq\int_DT^{\ast}(f)|g|\,d\mu +
\int_DT^{\ast}(f\chi_D)|g|\,d\mu<\infty.
\end{multline}
Suppose now in addition that $f$ is non-negative. Then by
\eqref{eq:3.1} we can define for $k \in \N$,
$$
S_kf(z) =
\left(\int_Df\,d\mu\right)^{-1}\int_D\int_{D^c}K(x,y)f(y)\,d\mu
yf(x)\,d\mu x
$$
when $z \in D \in \mathcal{D}_k$, where we interpret $S_kf(z) = 0$
when $z \in D \in \mathcal{D}_k$ with $\int_Df\,d\mu= 0$.

Let $\nu$ be the finite Borel measure on $X$ such that
$$
\nu(B) = \int_Bf\,d\mu
$$
for Borel sets $B \subset X$. Let $\mathcal A_k$ be the
$\sigma$-algebra generated by $\mathcal{D}_k$. We shall check that
$(S_kf,\mathcal{A}_k)$ is a martingale (with respect $\nu$).

Let $D \in \mathcal{D}_k$ and let $D_1,D_2,\dotsc$ be the sets in
$\mathcal{D}_{k+1}$ which form the disjoint partition of $D$. Then
\begin{equation*}
\begin{split}
\int_D S_{k+1}f\,d\nu &= \sum_i \int_{D_i}S_{k+1}f\,d\nu \\*[5pt] &=
\sum_i\int_{D_i}\frac{1}{\nu(D_i)}\int_{D_i}\int_{D_i^c}K(x,y)\,d\nu
y\,d\nu x \,d\nu\\*[5pt] &= \sum_i\int_{D_i}\int_{D_i^c}K(x,y)\,d\nu
y\,d\nu x\\*[5pt] &= \sum_i\int_{D_i}\sum_{j:i\neq
j}\int_{D_j}K(x,y)\,d\nu y\,d\nu x \\*[5pt] &\quad +
\sum_i\int_{D_i}\int_{D^c}K(x,y)\,d\nu y\,d\nu x\\*[5pt] &= 0 +
\int_D\int_{\R^n \setminus D}K(x,y)\,d\nu y\,d\nu x,
\end{split}
\end{equation*}
where $0$ comes from the antisymmetry of $K$.  This gives
$$
\frac{1}{\nu(D)}\int_DS_{k+1}f\,d\nu = S_kf(z) \text{ for } z \in D
$$
and implies that $(S_kf,\mathcal{A}_k)$ is a martingale.

Now we check that the martingale $(S_kf,\mathcal{A}_k)$ is
$L^1(\nu)$-bounded. We estimate using \eqref{eq:3.1}, the Schwartz
inequality and the $L^2$-boundedness of $T^{\ast}$,
\begin{equation*}
\begin{split}
\left|\int S_kf\,d\nu\right| \!&= \left|\sum_{D \in
\mathcal{D}_k}\frac{1}{\nu(D)}\int_D\int_{D^c}K(x,y)\,d\nu y\,d\nu
x\,\nu(D)\right|\\*[5pt] &=\left|\sum_{D \in
\mathcal{D}_k}\int_D\int_{D^c}K(x,y)f(y)\,d\mu yf(x)\,d\mu x
\right|\\*[5pt] &\leq \sum_{D \in
\mathcal{D}_k}\left(\int_DT^{\ast}(f)f\,d\mu +
\int_DT^{\ast}(f\chi_D)f\,d\mu\right)\\*[5pt] &\leq\sum_{D \in
\mathcal{D}_k}\left(\left(\int_D(T^{\ast}(f)^2\,d\mu\right)^{1/2}
\left(\int_Df^2\,d\mu\right)^{1/2}\right.\\*[5pt]
&\hskip3truecm\left.+
\left(\int_DT^{\ast}(f\chi_D)^2\,d\mu\right)^{1/2}\left(\int_Df^2\,d\mu\right)^{1/2}\right)\\*[5pt]
&\leq\!\sum_{D \in
\mathcal{D}_k}\!\!\left(\!\left(\int_D\!(T^{\ast}(f)^2\,d\mu\right)^{\!\!1/2\!}\!\!+\!
\left(C_0\int_D\!f^2\,d\mu\right)^{\!1/2}\right)\!\!\left(\int_D\!f^2\,d\mu\right)^{\!\!1/2\!}\\*[5pt]
&\leq\left(\left(\sum_{D \in
\mathcal{D}_k}\int_D(T^{\ast}(f)^2\,d\mu\right)^{1/2}
\right.\\*[5pt] &\hskip2truecm\left.+ \left(\sum_{D \in
\mathcal{D}_k}C_0\int_Df^2\,d\mu\right)^{1/2}\right) \left(\sum_{D
\in \mathcal{D}_k}\int_Df^2\,d\mu\right)^{1/2}\\*[5pt]
&=\left(\left(\int T^{\ast}(f)^2\,d\mu\right)^{\!1/2}\!+\!
\left(C_0\int f^2\,d\mu\right)^{\!1/2}\right) \left(\int
f^2\,d\mu\right)^{1/2}\\*[5pt] &\leq 2C_0^{1/2}\int f^2\,d\mu.
\end{split}
\end{equation*}
This proves the $L^1$-boundedness. Hence by the martingale
convergence theorem $(S_kf(z))$ converges for $\mu$ almost all $z\in
X$.

Now we assume also that
\begin{equation}\label{eq:3.4}
\lim_{k \to \infty}\sup\{\diam(D):D \in \mathcal{D}_k\} = 0.
\end{equation}
We define for $f \in L^2(\mu), k \in \N$,
$$
A_kf(z) = \frac{1}{\mu(D)}\int_D\int_{D^c}K(x,y)f(y)\,d\mu y\,d\mu x
\text{ when } z \in D \in \mathcal{D}_k,
$$
where $A_kf(z) = 0$ if $\mu(D) = 0$. Using the convergence of
$(S_kf(z))$ we shall now verify that for $f\in L^2(\mu)$ there
exists the finite limit
\begin{equation}\label{eq:3.5}
Tf(z) = \lim_{k \to \infty}A_kf(z)
\end{equation}
for $\mu$ almost all $z \in X$. Clearly, we may assume that $f$ is
non-negative. Moreover, since $A_k(f)=A_k(f+1)-A_k(1)$, we may
assume that $f\geq1$. To prove \eqref{eq:3.5} for such an $f$, write
$f_D=\frac{1}{\mu(D)}\int_Df\,d\mu$ for $D \in \mathcal{D}_k$ with
$\mu(D)>0$. Then by \eqref{eq:3.1}, the Schwartz inequality and
\eqref{eq:1.5} we have for $z\in D$,
\begin{equation*}
\begin{split}
&|S_kf(z)-A_kf(z)|\\*[5pt]
&=\left|\left(\int_Df\,d\mu\right)^{-1}\int_D\int_{D^c}K(x,y)f(y)\,d\mu
y\,(f(x)-f_D)\,d\mu x\right|\\*[5pt] &\leq
\frac{1}{\mu(D)}\left(\int_DT^{\ast}(f)|f-f_D|\,d\mu +
\int_DT^{\ast}(f\chi_D)|f-f_D|\,d\mu\right)\\*[5pt]
&\leq\frac{1}{\mu(D)}\left(\left(\int_DT^{\ast}(f)^2\,d\mu\right)^{1/2}\right.\\*[5pt]
&\hskip3.25truecm\left. +
\left(\int_DT^{\ast}(f\chi_D)^2\,d\mu\right)^{1/2}\right)\left(\int_D(f-f_D)^2\,d\mu\right)^{1/2}\\*[5pt]
&\leq \left(\frac{1}{\mu(D)}\left(2
\int_D(T^{\ast}(f)^2+C_0f^2)\,d\mu\right)^{1/2}
\left(\frac{1}{\mu(D)}  \int_D(f-f_D)^2\,d\mu\right)^{1/2}\right).
\end{split}
\end{equation*}
Here for $\mu$ almost all $z\in X$ as $k\to\infty$, the first factor
goes to $ 2^{1/2}\,(T^{\ast}(f)(z)^2+C_0f(z)^2)^{1/2}$, and the
second goes to $0$. Hence $S_kf(z)-A_kf(z)\to0$, which proves
\eqref{eq:3.5} for non-negative functions $f\in L^2(\mu)$ and of
course then also for all $f\in L^2(\mu)$. Moreover, $T\colon
L^2(\mu) \rightarrow L^2(\mu)$ is bounded.

To get from this the average convergence with balls one needs to
approximate balls with nested systems. At least in $\R^n$ this
approximation procedure can be done with dyadic cubes. The argument
is quite technical and will be omitted.

\section{Proof of Theorem 1.6}

We shall first make two reductions using the weak type
inequality~\eqref{eq:1.10}. Firstly, we may assume that $f=1$. To
see this, note that we may of course assume that $f$ is
non-negative. Bounded functions $f$ such that $f>\delta$ for some
$\delta>0$ are dense in the space of non-negative
$L^1(\mu)$-functions, whence standard techniques (as for
\eqref{eq:4.1} below) allow us to assume that $f$ is such a
function. Replacing $\mu$ by $f\mu$ gives then the reduction to
$f=1$.

Secondly, we may assume the uniform condition
\begin{equation}\label{eq:4.1}
\mu(B(x,r))\leq\eta(r)h(r)\leq h(r) \text{ for } x\in X,\quad r>0,
\end{equation}
where $\eta$ is a non-decreasing function such that $\eta(r)\to0$ as
$r\to0$. To see this, we use Egoroff's theorem  to select closed
sets $E_k,k=1,2,\dotsc$, such that $\mu(X\setminus E_k)<1/k$ and
$\mu(B(x,r))/h(r)\to0$ as $r\to0$ uniformly on $E_k$. Then using
\eqref{eq:1.10} we have for all $t>0$,
\begin{equation*}
\begin{split}
&\quad\;\mu(\{x:
\limsup_{\epsilon,\delta\to0}|T_{\epsilon}(1)(x)-T_{\delta}(1)(x)|>t\})\\*[3pt]
&=\mu(\{x:\limsup_{\epsilon,\delta\to0}
|T_{\epsilon}(1-\chi_{E_k})(x)-T_{\delta}(1-\chi_{E_k}(x)|>t\})\\*[3pt]
&\leq\mu(\{x:T^{\ast}(1-\chi_{E_k})(x)>t/2\})\\*[3pt] &\leq
C\mu(X\setminus E_k)/t,
\end{split}
\end{equation*}
provided the limit $\lim_{\epsilon\to0}T_{\epsilon}(\chi_{E_k})(x)$
exists for $\mu$ almost all $x\in E_k$. (It exists also for all
$x\in E_k^c$ since $E_k$ is closed.) That is, if we have the
convergence for the measures $\chi_{E_k}\mu$, which satisfy
\eqref{eq:4.1}, we have it also for $\mu$. Then it is easy to check
that the limit must be $T(1)(x)$ $\mu$ almost everywhere.

Thus it is enough to prove that
$\lim_{\epsilon\to0}T_{\epsilon}(1)(a)=T(1)(a)$ for $\mu$ almost all
$a\in X$ assuming \eqref{eq:4.1}. It is enough to consider points
$a\in X$ such that
$$
T(1)(a) =
\lim_{\epsilon\to0}\frac{1}{\mu(B(a,\epsilon))}\int_{B(a,\epsilon)}T(1)\,d\mu.
$$
Let $0<\delta<1/2$ and choose $p > 1/\delta$. Using \eqref{eq:2.2}
we can write for $\epsilon>0$,
\begin{equation*}
\begin{split}
\phi(\epsilon)&:= T_{\epsilon}(1)(a) -
\frac{1}{\mu(B(a,\epsilon))}\int_{B(a,\epsilon)}T(1)\,d\mu\\*[5pt]
&=\int_{B(a,p\epsilon)\setminus B(a,\epsilon)}K(a,x)\,d\mu x
\\*[5pt] &\quad+ \int_{B(a,p\epsilon)^c}K(a,x)\,d\mu x
+\frac{1}{\mu(B(a,\epsilon))}\int_{B(a,p\epsilon)^c}
T(\chi_{B(a,\epsilon)}) \,d\mu\\*[5pt]
&\quad+\frac{1}{\mu(B(a,\epsilon))}\int_{B(a,p\epsilon)\setminus
B(a,\epsilon)} T(\chi_{B(a,\epsilon)}) \,d\mu\\*[5pt]
&=\phi_1(\epsilon) + \phi_2(\epsilon) + \phi_3(\epsilon),
\end{split}
\end{equation*}
where
\begin{align*}
\phi_1(\epsilon)&=\int_{B(a,p\epsilon)\setminus
B(a,\epsilon)}K(a,x)\,d\mu x , \\*[5pt] \phi_2(\epsilon) &=
\int_{B(a,p\epsilon)^c}K(a,x)\,d\mu x +
\frac{1}{\mu(B(a,\epsilon))}\int_{B(a,p\epsilon)^c}
T(\chi_{B(a,\epsilon)}) \,d\mu  \intertext{and} \phi_3(\epsilon) &=
\frac{1}{\mu(B(a,\epsilon))}\int_{B(a,p\epsilon)\setminus
B(a,\epsilon)} T(\chi_{B(a,\epsilon)}) \,d\mu.
\end{align*}
The first term $\phi_1$ is easy to estimate:
\begin{equation*}
\begin{split}
|\phi_1(\epsilon)| &= \left|\int_{B(a,p\epsilon)\setminus
B(a,\epsilon)}K(a,x)\,d\mu x \right|\\*[5pt] &\leq
\frac{\mu(B(a,p\epsilon))}{h(\epsilon)}\leq
C_p\frac{\mu(B(a,p\epsilon))}{h(p\epsilon)} <\delta
\end{split}
\end{equation*}
by \eqref{eq:1.7} and \eqref{eq:1.9} for sufficiently small
$\epsilon$. Here and later $C_q$ for $q>1$ denotes a constant such
that $h(qr)\leq C_qh(r)$ for $r>0$. We estimate $\phi_2$ using
\eqref{eq:1.8} and \eqref{eq:4.1},
\begin{equation*}
\begin{split}
&|\phi_2(\epsilon)| \\*[5pt] &=
\left|\int_{B(a,p\epsilon)^c}K(a,x)\,d\mu x+
\frac{1}{\mu(B(a,\epsilon))}\int_{B(a,p\epsilon)^c}
T(\chi_{B(a,\epsilon)}) \,d\mu\right|\\*[5pt] &=
\left|\frac{1}{\mu(B(a,\epsilon))}\int_{B(a,\epsilon)}\!\left(\int_{B(a,p\epsilon)^c}\!K(a,x)\,d\mu
x\!-\!\int_{B(a,p\epsilon)^c}\!K(y,x)\,d\mu x\right)\,d\mu
y\right|\\*[5pt]
&\leq\frac{1}{\mu(B(a,\epsilon))}\int_{B(a,\epsilon)}\int_{B(a,p\epsilon)^c}
\frac{d(a,y)}{d(a,x)h(d(a,x))}\,d\mu x\,d\mu y\\*[5pt]
&\leq\epsilon\sum_{i=0}^{\infty}\frac{\mu(B(a,2^{i+1}p\epsilon))}{2^ip\epsilon
h(2^ip\epsilon)}\leq
\sum_{i=0}^{\infty}\frac{\mu(B(a,2^{i+1}p\epsilon))}{2^ipC_2^{-1}h(2^{i+1}p\epsilon)}
\leq 2C_2/p < 2C_2\delta.
\end{split}
\end{equation*}

To estimate $\phi_3$ we first show that at almost every point $\mu$
is doubling at some small scales. Then we only need to treat the
case $X=\R^n$. More precisely, let $C>2C_2$ be a constant and let
$F$ be the set of those~$a\in \R^n$ for which there exists
$\epsilon$, $0<\epsilon<1$, such that
$$
\mu(B(a,2^{1-k}\epsilon))\geq C\mu(B(a,2^{-k}\epsilon) \text{ for }
k=0,1,\dotsc.
$$
We also assume that $C>2^{n+1}$. We show now that $\mu(F)=0$. To
prove this we may assume that the support of $\mu$ is bounded, say
$\spt\mu\subset B(0,R)$. For $a\in F$ let $\epsilon=\epsilon(a)$ be
as above. Fix a big positive integer~$m$ and pick for each $a\in F$
an integer $k(a)\geq m$ such that for $k\geq k(a)$,
$$
C^{-k}\leq(2^{-k}\epsilon(a))^{n+1}.
$$
By Vitali's covering theorem (which holds in our setting, as we said
at the beginning of the introduction) we can find disjoint balls
$B(a_i,2^{-k_i}\epsilon_i)\subset B(0,R)$ with
$\epsilon_i=\epsilon(a_i)$ and $k_i\geq k(a_i)$ which cover $\mu$
almost all of $F$. Then
\begin{equation*}
\begin{split}
\mu(F)&\leq\sum_i\mu(B(a_i,2^{-k_i}\epsilon_i))\leq\sum_i
C^{-k_i}\mu(B(a,\epsilon_i))\\*[5pt]
&\leq\sum_i(2^{-k_i}\epsilon_i)^{n+1}\mu(\R^n) \leq
R^n2^{-m}\mu(\R^n).
\end{split}
\end{equation*}
Letting $m\to\infty$ we get $\mu(F)=0$.

Let now $a\in F^c$ and $0<\epsilon<1$. Then there is $k=0,1,\dotsc$,
such that
$$
\mu(B(a,2^{1-k}\epsilon))\leq C\mu(B(a,2^{-k}\epsilon))
$$
and
$$
\mu(B(a,2^{1-j}\epsilon))\geq C\mu(B(a,2^{-j}\epsilon)) \text{ for }
j=0,\dotsc,k-1,
$$
whence
$$
\mu(B(a,2^{-j}\epsilon))\leq C^{-j}\mu(B(a,\epsilon)) \text{ for }
j=0,\dotsc,k-1.
$$
Let $\epsilon_1=2^{-k}\epsilon$. Then $\mu(B(a,2\epsilon_1))\leq
C\mu(B(a,\epsilon_1))$ and, since $C>2C_2$, we get
\begin{equation*}
\begin{split}
|T_{\epsilon}(1)(a)-T_{\epsilon_1}(1)(a)|&\leq
\sum_{j=1}^k|T_{2^{1-j}\epsilon}(1)(a)-T_{2^{-j}\epsilon}(1)(a)|\\*[5pt]
&\leq\sum_{j=1}^k\int_{B(a,2^{1-j}\epsilon)\setminus
B(a,2^{-j}\epsilon)}|K(a,x)| \,d\mu x\\*[5pt]
&\leq\sum_{j=1}^k\frac{\mu(B(a,2^{1-j}\epsilon))}{h(2^{-j}\epsilon)}\\*[5pt]
&\leq\sum_{j=1}^k\frac{C^{1-j}\mu(B(a,\epsilon))}{C_2^{-j}h(\epsilon)}\leq
\,C \, \eta(\epsilon)\sum_{j=1}^k2^{-j}\leq
\,C\,\eta(\epsilon)<\delta
\end{split}
\end{equation*}
when $\epsilon$ is small enough. Consequently,
\begin{multline*}
|\phi(\epsilon)-\phi(\epsilon_1)|\leq
|T_{\epsilon}(1)(a)-T_{\epsilon_1}(1)(a)|\\*[5pt]
+\left|\frac{1}{\mu(B(a,\epsilon))}\int_{B(a,\epsilon)}T(1)\,d\mu-
\frac{1}{\mu(B(a,\epsilon_1))}\int_{B(a,\epsilon_1)}T(1)\,d\mu\right|<
\delta
\end{multline*}
when $\epsilon$ is small enough. Now we estimate the average of
$|\phi_3(t)|$ over $[\epsilon_1,2\epsilon_1]$ by
\begin{equation*}
\begin{split}
&\quad\;\frac{1}{\epsilon_1}\int_{\epsilon_1}^{2\epsilon_1}|\phi_3(t)|\,dt\\*[5pt]
&\leq
\frac{1}{\epsilon_1}\int_{\epsilon_1}^{2\epsilon_1}\frac{1}{\mu(B(a,t))}
\int_{B(a,pt)\setminus B(a,t)}\int_{B(a,t)}|K(x,y)|\,d\mu x\,d\mu
y\,dt \\*[5pt]
&=\frac{1}{\epsilon_1}\iiint_A\frac{1}{\mu(B(a,t))}|K(x,y)|\,d\mu
x\,d\mu y\,dt
\end{split}
\end{equation*}
where
\begin{multline*}
A=\{(x,y,t):d(x,a)<t\leq d(y,a)<pt,\, \epsilon_1\leq t\leq2\epsilon_1\}\\
\subset\{(x,y,t):d(x,a)<2\epsilon_1,\, d(y,a)<2p\epsilon_1,\,
d(x,a)<t\leq d(y,a)\}.
\end{multline*}
Thus by Fubini's theorem, \eqref{eq:1.7} and \eqref{eq:4.1},
\begin{equation*}
\begin{split}
\!&\!\frac{1}{\epsilon_1}\int_{\epsilon_1}^{2\epsilon_1}|\phi_3(t)|\,dt\\*[5pt]
\!&\!\leq\frac{1}{\epsilon_1\mu(B(a,\epsilon_1))}\!\int_{B(a,2p\epsilon_1)}\int_{B(a,2\epsilon_1)}|K(x,y)|
\int_{d(x,a)}^{d(y,a)}\,dt\,d\mu x\,d\mu y \\*[5pt]
\!&\!=\frac{1}{\epsilon_1\mu(B(a,\epsilon_1))}\!\int_{B(a,2\epsilon_1)}\int_{B(a,2p\epsilon_1)}
|K(x,y)|(d(y,a)-d(x,a))\,d\mu y\,d\mu x\\*[5pt] \!&\!
\leq\frac{1}{\epsilon_1\mu(B(a,\epsilon_1))}\!\int_{B(a,2\epsilon_1)}\int_{B(x,2(p+1)\epsilon_1)}|K(x,y)|\,d(x,y)\,d\mu
y\,d\mu x\\*[5pt] \!&\! \leq\frac{1}{\epsilon_1\mu
(B(a,\epsilon_1))}\!\!\int_{B(a,2\epsilon_1)}\!\sum_{i=0}^{\infty}\!\int_{\!B(x,2^{1-i}(p+1)\epsilon_1)\setminus
B(x,2^{-i}(p+1)\epsilon_1)}\!\!\! |K(x,y)|\,d(x,y)\,d\mu y\,d\mu x\!
\\*[5pt]
\!&\!\leq\frac{1}{\epsilon_1\mu(B(a,\epsilon_1))}\!\int_{B(a,2\epsilon_1)}\sum_{i=0}^{\infty}
\frac{2^{1-i}(p+1)\epsilon_1\mu(B(x,2^{1-i}(p+1)\epsilon_1))}{h(2^{-i}(p+1)\epsilon_1)}\,d\mu
x \\*[5pt]
\!&\!\leq\frac{1}{\epsilon_1\mu(B(a,\epsilon_1))}\!\sum_{i=0}^{\infty}
\frac{2^{1-i}(p\!+\!1)\epsilon_1\eta(2^{1-i}(p\!+\!1)\epsilon_1)h(2^{1-i}(p\!+\!1)\epsilon_1))}
{h(2^{-i}(p+1)\epsilon_1)} \mu(B(a,2\epsilon_1))\\*[5pt]
\!&\!\leq\frac{4C_2(p+1)\eta(2(p+1)\epsilon_1))\mu(B(a,2\epsilon_1))}{\mu(B(a,\epsilon_1))}
\leq 4CC_2(p+1)\eta(2(p+1)\epsilon)<\delta.
\end{split}
\end{equation*}
when $\epsilon$ is small enough. So there is $\epsilon_2,
\epsilon_1\leq\epsilon_2\leq2\epsilon_1$, such that
$|\phi_3(\epsilon_2)|<\delta$. Then
$|\phi(\epsilon_1)-\phi(\epsilon_2)|<\delta$ as above and so
\begin{multline*}
|\phi(\epsilon)|\leq|\phi(\epsilon)-\phi(\epsilon_1)|+|\phi(\epsilon_1)-\phi(\epsilon_2)|+|\phi(\epsilon_2)|\\
<2\delta+|\phi_1(\epsilon_2)|+|\phi_2(\epsilon_2)|+|\phi_3(\epsilon_2)|
<(4+2C_2)\delta.
\end{multline*}
This completes the proof of Theorem 1.6.

\section{Remarks on rectifiability}

One motivation for the developments in this paper was to find some
new insight to the following problem:

Let $m$ be an integer, $0<m<n$, and let $\mu$ be an $m$-dimensional
Ahlfors-David-regular Borel measure on $\R^n$, as in Section~1. For
$i=1,2,\dotsc,n$ let $T^{\ast}_i$ be the maximal operator related to
$\mu$ and the kernel $|x-y|^{-m-1}(x_i-y_i)$. Suppose that each
$T^{\ast}_i$ is bounded in $L^2(\mu)$. Does $\mu$ have to be
rectifiable, or even uniformly rectifiable in the sense of David and
Semmes?

By the rectifiability of $\mu$ we mean that there are
$m$-dimensional $C^1$\guio{surfaces} $M_1,M_2,\dotsc$ such that
$\mu(\R^n \setminus \cup_i M_i) = 0$. For the definitions of uniform
rectifiability, see~\cite{DS}.

If $m=1$, the answer to the above question is yes by~\cite{MMV}, and
the regularity assumptions on $\mu$ can be considerably relaxed, see
\cite{T2}. The problem is open for $m\geq2$.

It was shown in \cite{MPr}, see also \cite{M}, that the
rectifiability of an Ahlfors-David-regular measure~$\mu$ follows
from the existence of the principal values
$$
\lim_{\epsilon \to 0}\int_{\R^n \setminus
B(x,\epsilon)}|x-y|^{-m-1}(x_i-y_i)\,d\mu y,\quad  i=1,\dots,n,
$$
for $\mu$ almost all $x \in \R^n$. But it is not known if the
$L^2$-boundedness implies the above almost everywhere convergence.
Thus Theorem~\ref{thm:1.4} is a kind of replacement for this.
Unfortunately we don't know if the almost everywhere convergence of
the averages of Theorem~\ref{thm:1.4} implies rectifiablity, nor do
we know if it implies the almost everywhere existence of the
principal values in this particular case.

These questions are also related to geometric properties of
removable sets of bounded analytic functions in $\mathbb{C}$, see
\cite{MMV}, \cite{P} and \cite{T3}, and of Lipschitz harmonic
functions in $\R^n$, see \cite{MP}.

The $L^2$-boundedness does not always imply the almost everywhere
existence of principal values in the setting of Theorem~1.4. This
can be seen by considering a standard example of a purely
unrectifiable $1$\guio{dimensional} Ahlfors-David-regular set in the
plane which is the Cantor set obtained by starting with the unit
square, taking four squares of side-length~$1/4$ inside it in its
corners, then taking the squares of side-length~$1/16$ in the
corners of these, and so on. The final Cantor set~$C$ is the compact
set inside all these squares of all generations. In \cite{D2} David
constructed a $1$-dimensional odd Calder\'on-Zygmund kernel~$K$ such
that the operator $T^{\ast}$ related to~$K$ is bounded in $L^2(\mu)$
where $\mu$~is the natural ($1$-dimensional Hausdorff) measure
on~$C$.  However, it is easy to check that the principal values
$$
\lim_{\epsilon\to 0}\int_{B(x,\epsilon)^c}K(x-y)\,d\mu y
$$
fail to exist at $\mu$ almost all points $x \in \R^2$.

In \cite{H2} Huovinen considered homogeneous kernels such as
$$
K(z) = Re(z/|z|^2 - z^3/|z|^4)
$$
for $z\in \mathbb{C}$. He showed that there exist purely
unrectifiable $1$-dimensional  Ahlfors-David-regular sets on which
for such a kernel the principal values exist almost everywhere and
the related operator is bounded in $L^2$ on some subset of positive
measure. On the other hand, he showed in~\cite{H1} that for the
kernels $z^{2k-1}/|z|^{2k}$, $k=1,2,\dotsc$, and their linear
combinations the almost everywhere convergence of principal values
on $1$-dimensional AD-regular sets implies their rectifiability.

\begin{gracies}
 The authors are grateful to the ``Centre de Recerca Matem\`{a}tica" for the perfect organization of
  a special trimester in Fourier Analysis during the spring of 2006.
  The event offered an excellent opportunity to work on the
  paper.
 The second named author was partially supported by grants
 2005SGR00774 (Generalitat de Catalunya),  MTM2004-00519  and HF2004-0208.
\end{gracies}

\vspace{0.5 cm}

\begin{tabular}{l}

Pertti Mattila\\
Department of Mathematics and Statistics\\
P.O. Box 68 \\
FI-00014 \, University of Helsinki\\
Finland\\
 {\it E-mail:} {\tt Pertti.Mattila@Helsinki.fi}\\ \\
Joan Verdera\\
Departament de Matem\`{a}tiques\\
Universitat Aut\`{o}noma de Barcelona\\
08193 Bellaterra, Barcelona \\
 Catalonia\\
{\it E-mail:} {\tt jvm@mat.uab.cat}
\end{tabular}

\end{document}